\documentclass[english,]{article}
\usepackage{lmodern}
\usepackage{amssymb,amsmath}
\usepackage{ifxetex,ifluatex}
\usepackage{fixltx2e} 
\ifnum 0\ifxetex 1\fi\ifluatex 1\fi=0 
  \usepackage[T1]{fontenc}
  \usepackage[utf8]{inputenc}
\else 
  \ifxetex
    \usepackage{mathspec}
  \else
    \usepackage{fontspec}
  \fi
  \defaultfontfeatures{Ligatures=TeX,Scale=MatchLowercase}
\fi
\IfFileExists{upquote.sty}{\usepackage{upquote}}{}
\IfFileExists{microtype.sty}{%
\usepackage{microtype}
\UseMicrotypeSet[protrusion]{basicmath} 
}{}
\usepackage{hyperref}
\hypersetup{unicode=true,
            pdftitle={A Core Theory of Delay Systems},
            pdfauthor={Sébastien Boisgérault},
            pdfborder={0 0 0},
            breaklinks=true}
\ifnum 0\ifxetex 1\fi\ifluatex 1\fi=0 
  \usepackage[shorthands=off,main=english]{babel}
\else
  \usepackage{polyglossia}
  \setmainlanguage[]{english}
\fi
\usepackage{graphicx,grffile}
\makeatletter
\def\maxwidth{\ifdim\Gin@nat@width>\linewidth\linewidth\else\Gin@nat@width\fi}
\def\maxheight{\ifdim\Gin@nat@height>\textheight\textheight\else\Gin@nat@height\fi}
\makeatother
\setkeys{Gin}{width=\maxwidth,height=\maxheight,keepaspectratio}
\IfFileExists{parskip.sty}{%
\usepackage{parskip}
}{
\setlength{\parindent}{0pt}
\setlength{\parskip}{6pt plus 2pt minus 1pt}
}
\ifx\paragraph\undefined\else
\let\oldparagraph\paragraph
\renewcommand{\paragraph}[1]{\oldparagraph{#1}\mbox{}}
\fi
\ifx\subparagraph\undefined\else
\let\oldsubparagraph\subparagraph
\renewcommand{\subparagraph}[1]{\oldsubparagraph{#1}\mbox{}}
\fi

\title{A Core Theory of Delay Systems}
\author{Sébastien Boisgérault\footnote{MINES ParisTech, PSL Research University,
  Centre for Robotics, 60 Bd St Michel, 75006 Paris, France. E-mail:
  \texttt{Sebastien.Boisgerault@mines-paristech.fr}.}}
\date{}

\begin{document}
\maketitle

{
\setcounter{tocdepth}{3}
\tableofcontents
}
\newtheorem{theorem}{Theorem} \newtheorem{corollary}{Corollary}
\newtheorem{lemma}{Lemma} \newtheorem{definition}{Definition}

\newpage

\section*{Abstract}\label{abstract}
\addcontentsline{toc}{section}{Abstract}

We introduce a framework for the description of a large class of
delay-differential algebraic systems, in which we study three core
problems: first we characterize abstractly the well-posedness of the
initial-value problem, then we design a practical test for
well-posedness based on a graph-theoretic representation of the system;
finally, we provide a general stability criterion. We apply each of
these results to a structure that commonly arises in the control of
delay systems.

\newpage

\section{Introduction}\label{introduction}

Systems of delay-differential algebraic equations (DDAE) combine
differential and algebraic equations with delayed variables in the
right-hand side. Among them, we represent linear and time-invariant
systems with a finite memory \(r\) as

\begin{equation} \label{DDAE_}
  \begin{array}{rcl}
  \dot{x}(t) &=& A x_t + B y_t \\
        y(t) &=& C x_t + D y_t
  \end{array}
  \end{equation}

where \(x(t) \in \mathbb{R}^n\), \(y(t) \in \mathbb{R}^m\) and the
notation \(z_t\) stands for the memory of the variable \(z\) at time
\(t\): the function defined by \(z_t(\theta) = z(t+\theta)\) for any
\(\theta \in [-r,0]\). The symbols \(A\), \(B\), \(C\) and \(D\) denote
delay operators: linear and bounded operators from the space of
continuous functions \(C([-r,0], \mathbb{C}^{j})\) to \(\mathbb{C}^i\)
for some integers \(i\) and \(j\). Specifically, since
\((x(t), y(t)) \in \mathbb{R}^{n+m}\), we require that the compound
operator

\begin{equation} \label{compound}
  \left[
  \!
  \begin{array}{cc}
  A & B \\
  C & D 
  \end{array} 
  \!
  \right]
  \end{equation}

is linear and continuous from \(C([-r,0],\mathbb{C}^{n+m})\) to
\(\mathbb{C}^{n+m}\).

This description is general enough to encompass systems of seemingly
different nature -- such as ordinary differential equations (e.g.
\(\dot{x}(t) = x(t)\)), delay systems of retarded type (e.g.
\(\dot{x}(t) = x(t-1)\)) or neutral type (e.g.
\(\dot{y}(t) = \dot{y}(t-1)\), rewritten as \(\dot{x}(t) = 0\) and
\(y(t) = x(t) + y(t-1)\)), difference equations (e.g. \(y(t) = y(t-1)\))
but also different types of delays: systems with single delays (e.g.
\(\dot{x}(t) = x(t) + x(t-1)\)), multiple delays (e.g.
\(\dot{x}(t) = x(t-1) + x(t-2)\)), distributed delays (e.g.
\(\dot{x}(t) = \int_{t-1}^{t} x(\theta) \, d\theta\)) or a combination
thereof. We refer the reader to (Bensoussan et al.
\protect\hyperlink{ref-BPDM06}{2006}, chap. 4) for a more comprehensive
collection of examples.

However, competing theories of delay systems have been developped and
used over time (Hale \protect\hyperlink{ref-H77}{1977}; M. Delfour and
Karrakchou \protect\hyperlink{ref-DK87}{1987}; Salamon
\protect\hyperlink{ref-Sal84}{1984}); they may use different definitions
of delay operators, of the solutions to the initial-value problem and
may be restricted to systems of a certain type. For the newcomer in the
field, this is an unintended source of complexity. To make the subject
more widely accessible, we lay out in this paper a simple but general
framework for the description of delay systems, based on a combination
of linear algebra and measure theory. Then we develop on this foundation
a core theory: we characterize the well-posedness of the initial-value
problem, then we perform a graph-theoretic analysis of this issue and
finally, we provide a general stability criterion.

We illustrate the results of each section with the same example: the
so-called finite spectrum assignment (FSA) architecture which is used
for the control of dead-time systems (Manitius and Olbrot
\protect\hyperlink{ref-MO79}{1979}). This example -- and the methods
that we use to study it -- are actually representative of a large class
of controllers for delay systems which ranges from the primeval Smith
predictor to the general observer-predictor structure (Mirkin and Raskin
\protect\hyperlink{ref-MR03}{2003}). A linear dead-time systems is
governed by an ordinary differential equation
\(\dot{x}(t) = Ex(t) + F u(t)\) but the only information available at
time \(t\) is the delayed variable \(x(t-T)\). If we apply the control
\(u(t) = -G y(t)\) where \(y(t)\) is the predicted value of \(x(t)\)
based on \(x(t-T)\) and the history of \(u\) defined by \[
  y(t) = e^{T E} x(t-T) + \int_{-T}^0 e^{-\theta E} F  u_t(\theta) d\theta
  \] then the closed-loop dynamics is

\begin{equation} \label{FSA}
  \begin{array}{rll}
  \dot{x}(t) &=& E x(t) - FG y(t) \\
  y(t) &=& \displaystyle e^{T E} x(t-T) - \int_{-T}^0 e^{-\theta E} FG  y_t(\theta) d\theta
  \end{array}
  \end{equation}

Since this system is not overly simplistic (it combines effectively
differential and algebraic equations, its delays are discrete and
distributed), it is a good testbed for the theory exposed in this paper.

\section{The Initial Value Problem}\label{the-initial-value-problem}

\subsection{Delay Operators}\label{delay-operators}

We have described general delays as operators applied to continous
functions: if \(C(X,Y)\) denotes the space of continuous functions from
\(X\) to \(Y\), a delay operator with memory length \(r\) is a bounded
linear operator from \(C([-r,0],\mathbb{C}^j)\) to \(\mathbb{C}^i\) for
some integers \(i\) and \(j\). But delays also have two alternate (and
equivalent) representations -- as measures and as convolution kernels --
which we will use extensively in the next sections.

We start with a quick example of these representations: consider
\(y(t) = x(t-1)\) where \(x(t)\) is a scalar variable. We will associate
to this delayed expression three distinct representations \(m_1\),
\(m_2\) and \(m_3\) such that by construction

\begin{equation} \label{spec-ex}
  y(t) = m_1 x_t = 
  \int_{-1}^0 x_t(\theta) \, d m_2(\theta) 
  = \int_0^1 x(t - \theta) dm_3(\theta).
  \end{equation}

The first one \(m_1\) is an operator applied to continuous scalar
functions defined on \([-1,0]\); here, equation \eqref{spec-ex} clearly
mandates that \(m_1 \chi = \chi(-1)\). The representations \(m_2\) and
\(m_3\) are measures; equation \eqref{spec-ex} holds if \(m_2\) is
\(\delta_{-1}\), the Dirac measure at \(t=-1\) and \(m_3\) is
\(\delta_1\), the Dirac measure at \(t=1\). Now if instead we consider
the distributed delay \(y(t) = \int_{t-1}^t x(\theta) d\theta\),
equation \eqref{spec-ex} is satisfied with \[
  m_1 \chi = \int_{-1}^0 \chi(t) \, dt, \;
  m_2 = dt|_{[-1,0]} \; \mbox{ and } \;
  m_3 = dt|_{[0,1]}.
  \]

In general, we deal with vector-valued variables; in this broader
context, the term ``measure'' may refer to several things. We use the
term \emph{scalar measure} for a complex-valued and countably additive
function defined for the bounded Borel subsets of the real line (Fell
and Doran \protect\hyperlink{ref-FD88}{1988}). A \emph{matrix-valued
measure} (resp.~\emph{vector-valued measure}) is a countably additive
function defined for the bounded Borel subsets of the real line whose
values are complex matrices (resp.~vectors) of a fixed size.
Equivalently, it is as a matrix (resp.~vector) whose elements are scalar
measures.

Given a matrix-valued measure \(L\), we define the integral of the
vector function \(\phi:\mathbb{R} \to \mathbb{C}^j\) with respect to
\(L\) as the vector of \(\mathbb{C}^i\) whose \(k\)-th element is given
by \[
  \forall \, k \in \{1,\dots, j\}, \;
  \left[ \int dL \, \phi \right]_k
  = \sum_{l=1}^j \int\phi_l \, d L_{kl}.
  \] Riesz's representation theorem provides a bijection between delay
operators with memory length \(r\) and matrix-valued measures supported
on \([-r,0]\). We may therefore use the same symbol \(L\) to denote both
objects, and this convention yields for any
\(\phi \in C([-r,0], \mathbb{C}^j)\) \[
  L\phi = \int dL \, \phi.
  \] The nature of the argument of \(L\) (function or set) determines
without ambiguity which representation of \(L\) is used in a given
context. When the argument is a set, we will also drop the parentheses
when it improves readability: typically, instead of \(L(\{0\})\), we
will use the lighter notation \(L\{0\}\).

Let \(L^{\ast}\) refer to the measure obtained by symmetry of the
measure \(L\) around \(t=0\), which is defined for any bounded Borel set
\(B\) by \(L^{\ast}(B)= L(-B).\) This construct obviously provides a new
bijection between delay operators with memory length \(r\) and
matrix-valued measures supported on \([0,r]\), which is related to the
representation of delay operators as convolutions.

We say that a measure \(\mu\) is \emph{limited on the left} if its
support is included in \([-r, +\infty)\) for some \(r \in \mathbb{R}\),
\emph{causal} if its support is included in \([0, +\infty)\) and
\emph{strictly causal} if additionally \(\mu\{0\}=0\). The convolution
\(\mu \ast \nu\) of two scalar measures limited on the left \(\mu\) and
\(\nu\) is the scalar measure defined for any bounded Borel set \(B\) by
\[
  (\mu \ast \nu) (B) = \int \nu(B-x) d\mu(x).
  \] The convolution of two vector or matrix-valued measures limited on
the left combines the scalar convolution and the linear algebra product:
for example, the convolution of two measures \(L\) and \(M\) with values
in \(\mathbb{C}^{i\times j}\) and \(\mathbb{C}^{j\times k}\) is the
measure \(L \ast M\), with values in \(\mathbb{C}^{i \times k}\), given
by \[
   (L \ast M)(B)_{\alpha, \beta} = \sum_{\gamma =1}^{j} (L_{\alpha, \gamma} \ast M_{\gamma, \beta})(B).
   \] We may identify a (locally integrable) function \(f\) with the
measure given by \[
  f(B) = \int_B f(t) \, dt.
  \] This enables us to define the convolution of functions and
measures, which may be identified with a function, and of two functions,
which may be identified with an absolutely continous function.
Additionally, we implicitly extend by \(0\) functions that are defined
on a proper subset of real line -- for example functions defined on
\([-r,0]\) or \([-r, +\infty)\). With these conventions, if \(L\) is a
delay operator with memory \(r\) and
\(\phi \in C([-r,+\infty),\mathbb{C}^i)\) then for almost every \(t>0\)
\[
  L \phi_t  = (L^{\ast}  \ast \phi) (t).
  \]

\subsection{Well-Posedness}\label{well-posedness}

Several concepts of solution exist for the problem formally described by

\begin{equation} \label{DDAE}
  \forall \, t>0 \; \left|
  \begin{array}{rcl}
  \dot{x}(t) &=& A x_t + B y_t \\
        y(t) &=& C x_t + D y_t
  \end{array}
  \right. \! \!, \;
  (x(0), x_0, y_0) = (\phi, \chi, \psi).
  \end{equation}

A classic (or continuous) solution is a pair of continuous functions
\((x,y)\) from \([-r,+\infty)\) to \(\mathbb{C}^{n+m}\) such that
\(\dot{x}\) exists, is continous on \((0, +\infty)\) and satisfies the
system equations and initial conditions. In this setting, necessarily
the functions \(\chi\) and \(\psi\) are continous, \(\phi = \chi(0)\)
and the consistency condition \[
  \psi(0) = C \chi + D \psi
  \] holds.

The continuity assumption ensures that the right-hand side of the system
equations is defined for any time \(t\); if we relax this assumption, we
may generalize the concept of solution to locally integrable solutions
(Salamon \protect\hyperlink{ref-Sal84}{1984}; M. Delfour and Karrakchou
\protect\hyperlink{ref-DK87}{1987}).

Let \(L\) be a delay operator with memory length \(r\) and \(L^{\ast}\)
the corresponding convolution kernel. For any continuous function \(z\)
defined on \([-r,+\infty)\), we have \(L z_t = (L^{\ast} \ast z)(t)\)
for any \(t>0\) but the right-hand side of this equation is still
properly defined -- as a locally integrable function of \(t\) -- if
\(z\) is merely locally integrable, a strong incentive to rewrite the
initial value problem as a convolution equation. Using the Heaviside
step function \(e\), defined by \[
  e(t) = \left|
  \begin{array}{rl}
  1 & \mbox{if } \, t\geq 0, \\
  0 & \mbox{otherwise,}
  \end{array}
  \right.
  \] we may also rewrite the (integral form of) the differential
equation as a convolution equation. We end up with the following
definition:

\textbf{Definition -- Solution of the Initial Value Problem.} \emph{A
pair of locally integrable functions \((x, y)\), defined on
\([-r,+\infty)\), with values in \(\mathbb{C}^{n+m}\), is a (locally
integrable) solution of the DDAE initial value problem (\ref{DDAE}) if
\[
  (x(0^+), x_0, y_0) = (\phi, \chi, \psi)
  \] and if there is a \(f \in \mathbb{C}^n\) such that} \[
  \left[
  \begin{array}{cc}
  x \\ y
  \end{array}
  \right](t)
  = 
  \left[
  \begin{array}{cc}
  e \ast A^{\ast} & e \ast B^{\ast} \\
  C^{\ast} & D^{\ast}
  \end{array}
  \right]
  \ast 
  \left[
  \begin{array}{c}
  x \\ y
  \end{array}
  \right](t)
  +
  \left[
  \begin{array}{c}
  f \\ 0
  \end{array}
  \right]
  \; \; 
  \mbox{\it for a.e.} \; t>0.
  \]

This definition makes sense even if the initial function data is not
continuous, does not meet the consistency condition, or if
\(\chi(0)\neq \phi\). It is however consistent with the concept of
continuous solution when such solutions exist. The constant vector \(f\)
is uniquely determined by the initial value: we have \(x(0^+) =\phi\) if
and only if
\(f = \phi - (e\ast A^{\ast} \ast \chi + e\ast B^{\ast} \ast \psi)(0)\).

Now, because of its algebraic component, DDAE system \eqref{DDAE} may
have no solutions or multiple solutions; we may actually easily exhibit
an algebro-differential system (with no delay) with this property.
Consider

\begin{equation} \label{ADE}
  \forall \, t > 0, \;
  \left|
  \begin{array}{ccl}
  \dot{x}(t) &=& 0 \\
  y(t) &=& x(t) + y(t)
  \end{array}
  \right.
  \end{equation}

It is defined formally by \(n=1\), \(m=1\), for example \(r=1\) (any
nonnegative value is admissible) and for any
\(\chi \in C^0([-r,0],\mathbb{R}^{n})\) and
\(\psi \in C^0([-r,0],\mathbb{R}^{m})\), \(A\chi = 0\), \(B\psi=0\),
\(C\chi=\chi(0)\), and \(D\psi = \psi(0)\). If \(x(0^+) = 0\), then
\(x(t) = 0\) for \(t>0\) is the unique solution of the first system
equation; the second equation becomes \(y(t) = y(t)\), hence arbitrary
values of \(y(t)\) for \(t>0\) satisfy it: multiple solutions exist. On
the contrary, if \(x(0^+)=1\), then necessarily \(x(t) = 1\) for any
\(t>0\); the second equation would become \(y(t) = 1 + y(t)\), which no
function \(y(t)\) may satisfy: there are no solutions to this initial
value problem.

A simple assumption that yields uniqueness of the solution is
\emph{explicitness}: a DDAE system is explicit if the right-hand side of
its algebraic equation only depends strictly causally on the variable
\(y\), that is, if \(D\{0\} = 0\). For example, system \eqref{ADE} is
not explicit since \(y(t) = x(t) + y(t)\) and thus \(D\{0\} = 1\);
replace this equation with \(y(t) = x(t) + y(t-\epsilon)\) for any
\(\epsilon>0\) and it becomes explicit. However this assumption is too
conservative for some practical use cases, including the kind of sound
composition of input-output systems that is exposed in the next section.
Therefore we state in this section a well-posedness theorem applicable
under a more general assumption that ensures that the system is merely
equivalent to an explicit system.

Let \(I_p\) be the \(p \times p\) identity matrix and \(X\) be the
product space

\begin{equation}
  X= \mathbb{C}^{n} \times L^2([-r,0],\mathbb{C}^{n+m}).
  \end{equation}

endowed with the norm

\begin{equation} \label{normX}
  \|(\phi, \chi, \psi)\|_X^2
  =
  |\phi|^2 + \int_{-r}^0 |\chi(t)|^2 + |\psi(t)|^2 \, dt.
  \end{equation}

\textbf{Theorem -- Well-posedness.} \emph{DDAE systems such that
\(I_m - D\{0\}\) is invertible are well-posed in the product space
\(X\): when \((\phi, \chi, \psi) \in X\), there is a unique solution
\((x,y)\) to the initial value problem \eqref{DDAE}, it belongs to
\(L^2_{\rm loc}([-r,+\infty),\mathbb{C}^{n+m})\), \(x\) is continuous
and for any \(t>0\) there is some \(\alpha>0\) such that}

\begin{equation} \label{WPI}
  \|(x(t), x_t, y_t)\|_X
  \leq 
  \alpha
  \|(\phi, \chi, \psi)\|_X.
  \end{equation}

It is classic that this result holds when the system is explicit
(Salamon \protect\hyperlink{ref-Sal84}{1984}; M. Delfour and Karrakchou
\protect\hyperlink{ref-DK87}{1987}; Boisgérault
\protect\hyperlink{ref-Boi13}{2013}). Therefore, to prove the general
case, we only need to demonstrate the following lemma:

\textbf{Lemma -- Equivalent Systems.} \emph{Assume that the matrix
\(J := I_m - D\{0\}\) associated to system \eqref{DDAE} is invertible
and denote \(E\) and \(F\) the delay operators \(E := J^{-1} C\) and
\(F := J^{-1} (D - D|_{\{0\}})\). The original initial value problem and
the one defined by}

\begin{equation} \label{DDAEX}
  \forall \, t>0 \; \left|
  \begin{array}{rcl}
  \dot{x}(t) &=& A x_t + B y_t \\
        y(t) &=& E x_t + F y_t
  \end{array}
  \right. \! \!, \;
  (x(0), x_0, y_0) = (\phi, \chi, \psi).
  \end{equation}

\emph{have the same solutions but the latter system is always explicit:}
\[
  F\{0\} = 0.
  \]

\textbf{Proof.} We may decompose the convolution kernel of \(D\) into
two components: \[
  D^{\ast} = (D|_{[-r,0]})^{\ast} = (D|_{\{0\}})^{\ast} + (D|_{[-r, 0)})^{\ast}.
  \] The first component is instantaneous: for any locally integrable
function \(y\) \[
  ((D|_{\{0\}})^{\ast} \ast y)(t) = D\{0\} y(t);
  \] the second component is strictly causal: \[
  (D|_{[-r, 0)})^{\ast} \{0\} = D|_{[-r, 0)} \{0\}= 0.
  \] Consequently, a pair of locally integrable functions
\((x,y): [-r,+\infty) \to \mathbb{R}^{n+m}\) is a solution of the DDAE
system equations \eqref{DDAE} if and only if they satisfy, almost
everywhere for \(t>0\), the equations \[
  \begin{array}{rcl}
  x(t) &=& (A^{\ast} \ast x)(t) + (B^{\ast} \ast y)(t) + f \\
  y(t) - D\{0\} y(t) &=& (C^{\ast} \ast x)(t) + ((D|_{[-r, 0)})^{\ast} \ast y)(t)
  \end{array}
  \] By assumption the matrix \(J = I_m - D\{0\}\) is invertible, thus
the second equation is equivalent to
\[y(t) = (E^{\ast} \ast x)(t) + (F^{\ast} \ast  y)(t)\] with
\(E = J^{-1} C\) and \(F = J^{-1} D|_{[-r, 0)}\). The functions \(x\)
and \(y\) are solutions of the original DDAE system if and only if they
are solutions -- with the same initial values -- of the system whose
delay operators are \(A\), \(B\), \(E\) and \(F\). This new DDAE system
is explicit by construction: \(F\{0\} = J^{-1} D|_{[-r, 0)} \{0\}=0\).
\hfill \(\blacksquare\)

\subsection{The FSA System}\label{the-fsa-system}

The delay operators \(A\), \(B\) and \(C\) of system \eqref{FSA} are
defined by \(Ax_t = E x(t)\), \(B y_t = - FG y(t)\) and
\(C x_t = e^{TE} x(t-T)\). Since the Dirac measure \(\delta_{\tau}\) at
\(\tau \in \mathbb{R}\) satisfies
\(\delta_{\tau} \varphi = \varphi(\tau)\), we have \(A = E \delta_0\),
\(B = -FG \delta_0\), \(C = e^{TE}\delta_{-T}\) or equivalently

\begin{equation} \label{FSA-ker-1}
  A^{\ast} = E \delta_0, \; 
  B^{\ast} = -FG \delta_0, \; 
  C^{\ast} = e^{TE}\delta_{T}.
  \end{equation}

The fourth delay operator \(D\) is defined by \[
  D y_t 
  = - \int_{0}^T e^{\theta E} FG  y(t-\theta) d\theta
  = - (e^{\theta E} FG d\theta|_{[0,T]} \ast y)(t),
  \] hence

\begin{equation} \label{FSA-ker-2}
  D^{\ast} = - e^{\theta E} FG d\theta|_{[0,T]}.
  \end{equation}

In particular \[
  D^{\ast}\{0\}
  =
  - \int_{\{0\}} 
  e^{\theta E} FG d\theta = 0,
  \] thus this DDAE system is explicit and therefore it is well-posed.

\section{\texorpdfstring{Graph-Theoretic Analysis
\label{CC}}{Graph-Theoretic Analysis }}\label{graph-theoretic-analysis}

\begin{figure}[htbp]
\centering
\includegraphics{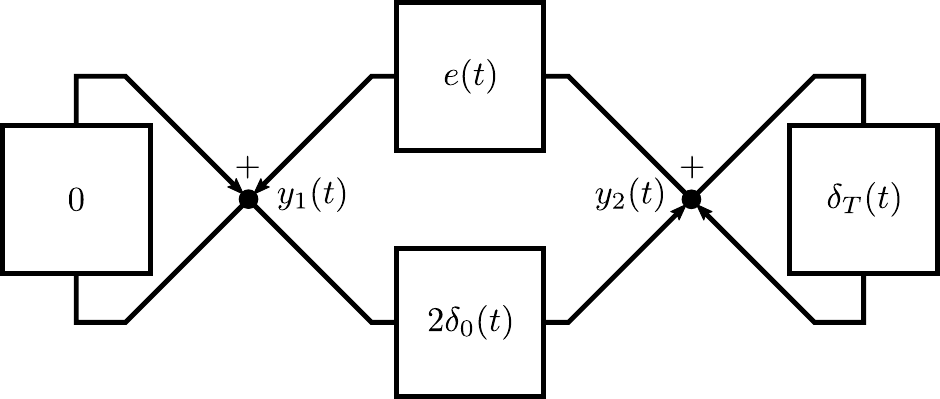}
\caption{Block diagram of a simple delay system. \label{sg}}
\end{figure}

Block diagrams used in control theory provide an alternative to the
system-of-equations approach for the modelling of dynamical systems: in
this context, dynamical systems are specified as diagrams made of blocks
and wires as in figure \ref{sg}. Blocks are independent input-output
systems with their own behaviors: each block specifies how to compute
the values of its output variables given the values of its input
variables (and optionally of some local initial state). In figure
\ref{sg}, each rectangle is a block and the symbol it holds defines its
behavior. Wires connect these blocks to build larger systems: they
specify how the input of each block is computed as a combination of the
outputs of all blocks.

This classic description -- actually a graph representation in disguise
-- offers new indsights into the structural properties of dynamical
systems. In section \ref{BD-sec}, we complete for delay systems the
informal definition of block diagrams that we have given so far and
expose the connexion with graph theory. Then, in section \ref{CA-sec},
we focus on causality to design sound composition rules of subsystems,
that ensure the well-posedness of the initial-value problem for the
global system.

\subsection{\texorpdfstring{Block Diagrams
\label{BD-sec}}{Block Diagrams }}\label{block-diagrams}

For linear time-invariant systems, block diagrams can always be
decomposed into elementary input-output systems with a scalar input
\(u(t)\) and a scalar output \(y(t)\), defined for \(t>0\), and
characterized by a scalar convolution kernel \(k^{\ast}\). Additionally,
for delay systems, only two types of elementary subsystems are required:
the integrator, whose kernel is the Heaviside step function \(e\), and
(scalar) finite-memory delay operators. By definition, the behavior of
an elementary system is governed by \[
  y(t) = (k^{\ast} \ast u)(t) + (k^{\ast} \ast \omega)(t), \; t > 0;
  \] where \(\omega\), a scalar function defined on \([-r,0]\) that
represents a finite memory of past input values, determines the system
initial state.

A block diagram embeds a collection of such systems into a (directed,
labeled) graph structure: each vertice \(i\) refers to an output
variable \(y_i\) and each edge \(j \to i\) is labeled with a scalar
kernel \(k^{\ast} = K^{\ast}_{ij}\); the global initial state is
determined by the collection of functions \(\omega_{ij}\) for every edge
\(j \to i\). The functional analytic representation of this system is

\begin{equation}
  y_i(t) = \sum_{j}(K_{ij}^{\ast} \ast y_j)(t) 
           +  \sum_{j}(K_{ij}^{\ast} \ast \omega_{ij})(t), \; t>0.
  \end{equation}

For example, the block diagram depicted in figure \ref{sg} has two nodes
-- associated to variables \(y_1\) and \(y_2\) -- and four edges --
whose labels are \(K^{\ast}_{11} = 0\), \(K^{\ast}_{12} = e\),
\(K^{\ast}_{21} = 2 \delta_0\) and \(K^{\ast}_{22} = \delta_T\). Then
\(y_1(t)\), defined for \(t>0\), satisfies \[
  \begin{split}
  y_1(t) &= (0 \ast y_1)(t) + (e \ast y_2)(t) + (0 \ast \omega_{11})(t) + (e \ast \omega_{12})(t) \\
         &= (e \ast y_2)(t) + (e \ast \omega_{12})(t)
  \end{split}
  \] which, since \((e\ast \omega_{12})(t)\) is a constant \(c\) for
\(t>0\), is the integral form of the differential equation
\(\dot{y}_1(t) = y_2(t)\) with initial value \(y_1(0) = c\). On the
other hand \(y_2(t)\), also defined for \(t>0\), satisfies \[
  y_2(t) = (2\delta_0 \ast y_1)(t) + (\delta_T \ast y_2)(t) +
(2\delta_0 \ast \omega_{21})(t) + (\delta_T \ast \omega_{22})(t)
  \] and since
\((2\delta_0 \ast \omega_{21})(t) = (\delta_T \ast \omega_{22})(t) = 0\)
for any \(t > r\), eventually the equation
\(y_2(t) = 2y_1(t) + y_2(t-T)\) holds and the behavior of the system is
ruled by a DDAE system. This is not accidental, since the following
result holds:

\textbf{Theorem.} \emph{Every DDAE system has a canonical block diagram
representation.}

\textbf{Proof.} Suppose that for some constant vector \(f\) \[
  \left|
  \begin{array}{l}
  x(t) = (e \ast A^{\ast} \ast x)(t) + (e \ast B^{\ast} \ast y)(t) + f \\
  y(t) = (C^{\ast} \ast x)(t) + (D^{\ast} \ast y)(t)
  \end{array}
  \right., \; t > 0
  \] and \[
  (x(0^+), x|_{[-r,0]}, y|_{[-r, 0]}) = (\phi, \chi, \psi). 
  \] The vector \(f\) is uniquely determined by the initial state. We
may define the auxiliary variable \(w(t)\) for \(t>0\) by
\(w(t) = (A^{\ast} \ast x)(t) + (B^{\ast} \ast y)(t),\) and rewrite the
first equation of the DDAE system as \[
  x(t) = (e \ast w)(t) + (e \ast \upsilon)(t), \; t>0
  \] for any function \(\upsilon:[-r,0] \mapsto \mathbb{C}^n\) such that
\((e \ast \upsilon)(0)=f.\) It is now plain that the variable
\(z(t) = (x(t), w(t), y(t))\) defined for \(t>0\), the matrix kernel and
the initial state

\begin{equation} \label{BD-DDAE1}
  K^{\ast}
  =
  \left[
  \begin{array}{ccc}
  0 & e I_n & 0 \\ 
  A^{\ast} & 0 & B^{\ast} \\
  C^{\ast} & 0 & D^{\ast} \\
  \end{array}
  \right], \;
  \omega_{ij} = (\chi, \upsilon, \psi)_j
  \end{equation}

satisfy

\begin{equation} \label{BD-DDAE3}
  z_i(t) = \sum_{j}(K_{ij}^{\ast} \ast z_j)(t) 
           +  \sum_{j}(K_{ij}^{\ast} \ast \omega_{ij}), \; t>0.
  \end{equation}

which is the functional analytic representation of a block diagram.
\hfill \(\blacksquare\)

\subsection{\texorpdfstring{Causality Analysis
\label{CA-sec}}{Causality Analysis }}\label{causality-analysis}

We say that a scalar input-output system is \emph{causal}
(resp.~\emph{strictly causal}) if its convolution kernel is causal
(resp.~strictly causal). For example, since the kernel of the integrator
is the function \(e\), it is strictly causal: \[
  e\{0\} = \int_{\{0\}} e(t)\, dt = 0.
  \] A (discrete or pure) delay of \(T>0\) seconds satisfies
\(y(t) = u(t-T)\); its kernel is \(k^{\ast} = \delta_{T}\) and it is
strictly causal. Similarly, a distributed delay with bounded delay
\(r\), defined by the equation \[
  y(t) = \int_{-r}^0 h(\theta) u(t+\theta) d\theta
  \] for some locally integrable function \(h:[-r, 0] \to \mathbb{R}\),
is also strictly causal. A gain \(y(t) = k u(t)\) is causal, but it is
strictly causal only if \(k=0\).

\textbf{Definition -- Causality Loop.} \emph{A causality loop of a block
diagram with kernel matrix \(K^{\ast}\) of size \(n\times n\) is a
finite sequence of integers \(i_0, i_1, \dots, i_j\) in
\(\{1, \dots, n\}\) such that \(i_0 = i_j\) and for any
\(\ell \in \{0,\dots, j-1\}\), the system with kernel
\(k^{\ast}_{\ell} = K^{\ast}_{i_{\ell+1}i_{\ell}}\) is not strictly
causal.}

The block diagram of figure \ref{sg} has no causality loop, since the
only element of its matrix \(K^{\ast}\) which is not strictly causal is
\(K^{\ast}_{21} = 2\delta_0\) which is not on the matrix diagonal. The
following result is applicable in this case:

\textbf{Theorem.} \emph{If the block diagram representation
\eqref{BD-DDAE1}-\eqref{BD-DDAE3} of DDAE system \eqref{DDAE} has no
causality loop then \(D\{0\}\) is nilpotent.}

\textbf{Proof.} Let \(K^*\) be the kernel matrix of the block diagram
representation of DDAE system \eqref{DDAE}. In this proof, for any
matrix-valued measure \(M\), we denote \(M_{0}\) the matrix \(M\{0\}\).
Let \(\mathcal{A}[K]\) be the adjacency matrix of the block diagram
directed graph, stripped off its strictly causal edges: \[
  \mathcal{A}[K]_{ij} = 
  \left|
  \begin{array}{ll}
  1 & \mbox{if } \, [K_0]_{ij} \neq 0 \\
  0 & \mbox{otherwise.}
  \end{array}
  \right.
  \] The graph-theoretic concept of loop corresponds to our definition
of causality loop. Element \((i,j)\) of the \(p\)-th power
\(\mathcal{A}[K]^p\) of the adjacency matrix is the number of distinct
paths with \(p\) edges between vertices \(j\) and \(i\). If there is a
loop of length \(p\) in the block diagram, there is also a loop whose
length is an arbitrary multiple of \(p\), hence the adjacency matrix is
not nilpotent. Conversely, if the adjacency matrix is not nilpotent,
there is a pair \((i,j)\) such that \(\mathcal{A}[K]^n_{ij} \neq 0\);
since any path between \(j\) and \(i\) goes through \(n+1\) vertices and
the graph has only \(n\) distinct vertices, at least one of them is
repeated: the path necessarily contains a loop.

Now, if \(\mathcal{A}[K]^p = 0\) for some \(p \in \mathbb{N}^*\) then
\([K_0]^p=0\). Indeed, since \[
  \mathcal{A}[K]^p_{ij} = \sum_{\ell_1, \dots, \ell_{p-1}}  
  \mathcal{A}[K]_{i\ell_1} \dots \mathcal{A}[K]_{\ell_{p-1}j}
  \] and all elements of the adjacency matrix are non-negative, if its
\(p\)-th power is zero, there is for each term in the sum above at least
one factor \(\mathcal{A}[K]_{\nu\mu}\) which is zero and thus
\([K_0]^p_{\nu\mu} = 0\). Consequently, since \[
  [K_0]^p_{ij} = \sum_{\ell_1, \dots, \ell_{p-1}}  
  [K_0]_{i\ell_1} \dots [K_0]_{\ell_{p-1}j},
  \] the \(p\)-th power of \(K_0\) is also zero. But for any integer
\(p \geq 2\), we have \[
  [K_{0}]^p
  =
  \left[
  \begin{array}{ccc}
  0 & 0 & 0     \\
  A_0 & 0 & B_0 \\
  C_0 & 0 & D_0 \\
  \end{array}
  \right]^p
  =
  \left[
  \begin{array}{ccc}
  0 & 0 & 0 \\
  B_0 [D_0]^{p-2} C_0 & 0 & B_0 [D_0]^{p-1} \\
  \![D_0]^{p-1} C_0 & 0 & [D_0]^{p} \\
  \end{array}
  \right].
  \] Consequently, if the block diagram has no causality loop,
\(\mathcal{A}[K]\) is nilpotent and the matrix \(D_0=D\{0\}\) is also
nilpotent. \hfill \(\blacksquare\)

Note that the absence of causality loop does not necessarily provide
explicitness: for example, the delay operator \(D\) of the DDAE system
governed by \(y_1(t) = y_2(t)\) and \(y_2(t) = y_1(t-1)\) satisfies \[
  D\{0\} = 
  \left[
  \begin{array}{cc}
  0 & 1 \\
  0 & 0
  \end{array} 
  \right]
  \] which is not zero. However, it is plain that its block-diagram
representation has no causality loop: it contains only one subsystem
which is not strictly causal and this subsystem connects variable
\(y_1\) to the other variable \(y_2\). This is sufficient to ensure its
well-posedness:

\textbf{Corollary.} \emph{A DDAE system without causality loop is
well-posed.}

\textbf{Proof.} If DDAE system \eqref{DDAE} has no causality loop, the
matrix \(D\{0\}\) is nilpotent and \([D\{0\}]^m = 0\). Consequently, the
matrix \(I_m - D\{0\}\) is invertible:
\[[I_m - D\{0\}]^{-1} = I_m + D\{0\} + \dots + [D\{0\}]^{m-1}.\] The
well-posedness theorem of the previous section provides the conclusion.
\hfill \(\blacksquare\)

\subsection{Block Diagram of the FSA
System}\label{block-diagram-of-the-fsa-system}

\begin{figure}[htbp]
\centering
\includegraphics{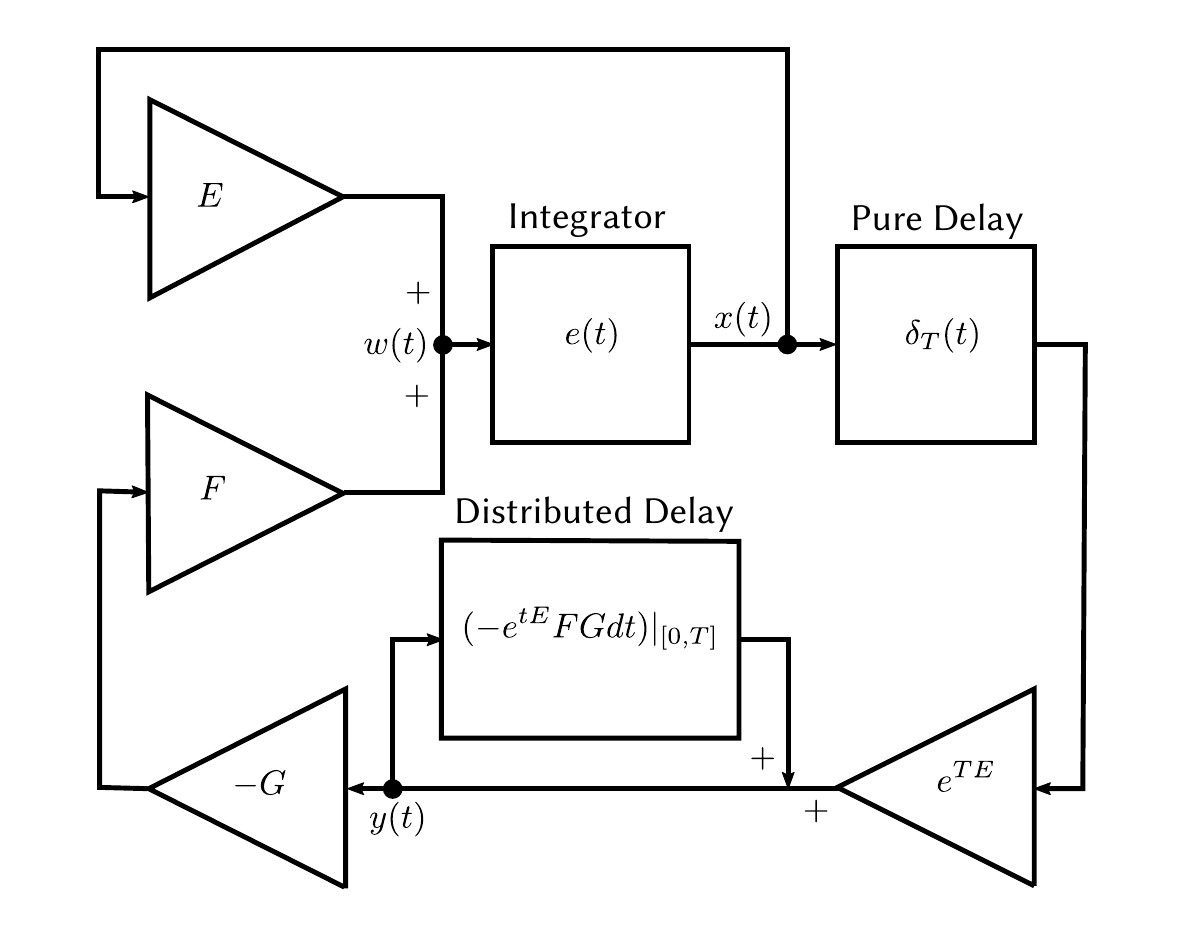}
\caption{FSA closed-loop system block diagram. Each rectangular block
characterizes an input-output system by its (function or measure)
convolution kernel. A triangle with symbol \(K\) represents a gain: a
system whose output \(y(t)\) and input \(u(t)\) satisfy
\(y(t) = K u(t)\). \label{FSAdiagram}}
\end{figure}

In the block diagram representation of system \eqref{FSA}, the auxiliary
variable \(w\) satisfies \[
  w(t) = E x(t) - FG y(t), \; t > 0.
  \]

Since \(x(t) = e \ast w(t) + e \ast \upsilon (t)\) is the only (system
of) equation(s) with some variable \(w_j\) in the right-hand side or
some variable \(x_i\) in the left hand-side, there is a unique outgoing
edge of \(w_j\) (with a nonzero label) which is also the unique incoming
edge of the variable \(x_j\). This label is the function \(e\), strictly
causal, thus no causality loop may contain either of the components of
the state variables \(w\) or \(x\). Now the possible causality loops
that could remain would contain only vertices which are components of
the state variable \(y\). But the edge between \(y_j\) and \(y_i\) is
\(H^{\ast}_{ij}\) where \(H^{\ast} = (-e^{tE} FG dt) |_{[0,T]}\). Again,
this kernel is a function and hence is strictly causal. Thus no such
loop may exist and the initial-value problem is well-posed.

This analysis may also be carried out graphically on figure
\ref{FSAdiagram}. Once the integrator, the pure delay and the
distributed delay -- which are strictly causal subsystems -- are removed
from the diagram, it is plain that there is no loop and thus that the
system has no causality loop.

\section{Stability}\label{stability}

\subsection{Test in the Laplace
Domain}\label{test-in-the-laplace-domain}

There is a general and simple test in the Laplace domain for the
exponential stability of a DDAE system. It can also be used to determine
its growth bound, defined as the infimum of the \(\sigma^+\) such that
for any initial value, there is a \(\kappa > 0\) such that the solution
satisfies \(\|(x(t), x_t, y_t)\|_X \leq \kappa e^{\sigma^+ t}\) for any
\(t>0\); exponential stability corresponds to a negative growth bound.

A few definitions are in order first: the Laplace transform
\(\mathcal{L}L\) of a compactly supported measure \(L\) is the function
defined on \(\mathbb{C}\) by \[
    \mathcal{L} L (s) = \int e^{-st} \, dL(t).
    \] The Laplace transform of \(L^{\ast}\) satisfies
\(\mathcal{L} L^{\ast} (s) = L (t \mapsto e^{st}).\) The
\emph{characteristic matrix} of the system is then defined as

\begin{equation} \label{characteristic}
  \Delta(s) =   
  \left[
    \begin{array}{cc}
    sI_n & 0 \\
     0 & I_m
    \end{array}
  \right] 
  -   
  \mathcal{L} 
  \left[
    \begin{array}{cc}
    A^{\ast} & B^{\ast} \\
    C^{\ast} & D^{\ast}
    \end{array}
    \right](s)
  \end{equation}

and its \emph{characteristic function} as its determinant.

We may now state the criterion:

\textbf{Theorem.} \emph{The growth bound \(\sigma\) of a DDAE system
without causality loop is:}

\begin{equation}
  \sigma
  = 
  \sup \; \{\Re (s) \; | \; \, s \in \mathbb{C}, \; \det \Delta(s) = 0 \}.
  \end{equation}

The proof of this theorem is carried out in (Boisgérault
\protect\hyperlink{ref-Boi13}{2013}) in the case of explicit systems.
Here, we expose the proof in the general case.

\textbf{Lemma.} \emph{DDAE systems \eqref{DDAE} and \eqref{DDAEX} have
simultaneously singular characteristic matrices \(\Delta(s)\) and
\(\Delta^{\dagger}(s)\).}

\textbf{Proof.} Let \(J = I_m - D\{0\}\), \(E = J^{-1} C\) and
\(F = J^{-1}(D - D|_{\{0\}})\); we have \[
  \mathcal{L} E^{\ast} = J^{-1}\mathcal{L} C^{\ast}
  \; \mbox{ and } \;
  \mathcal{L} F^{\ast} = J^{-1}(\mathcal{L} D^{\ast} - D\{0\})
  \] and since \(I_m + J^{-1} D \{0\} = J^{-1}\),
\(I_m - \mathcal{L} F^{\ast} = J^{-1}(I_m - \mathcal{L}D^{\ast})\).
Finally, \[
  \Delta^{\dagger}(s)
  =
    \left[
    \begin{array}{cc}
    sI_n & 0 \\
     0 & I_m
    \end{array}
  \right] 
  -   
  \mathcal{L} 
  \left[
    \begin{array}{cc}
    A^{\ast} & B^{\ast} \\
    E^{\ast} & F^{\ast}
    \end{array}
    \right](s)
  =
  \left[
    \begin{array}{cc}
    I_n & 0 \\
     0 & J^{-1}
    \end{array}
  \right] 
  \Delta(s)
  \] and the characteristic matrices are simultaneously singular.
\hfill \(\blacksquare\)

The determinant and adjugate of the characteristic matrix
\(\Delta^{\dagger}\) have a quasi-polynomial structure:

\begin{equation} \label{poly}
    \det \Delta^{\dagger} (s) = \sum_{i=0}^n c_i(s) s^i, \; \;
    \mathrm{adj}\, \Delta^{\dagger} (s) = \sum_{i=0}^{n} C_i(s) s^i
    \end{equation}

where the \(c_i\) (resp. \(C_i\)) are entire functions (resp. matrices
of entire functions) bounded on any right-hand plane. It is plain that
the leading coefficient of the characteristic function is given by
\(c_{n}(s) = \det \Delta_0(s)\) where \(\Delta_0\) is the characteristic
matrix of the system \(y(t) = F y_t\).

\textbf{Lemma -- Zero Clusters.} \emph{Let \(H_{\sigma}\) be the open
right half-plane \(\{s \in \mathbb{C} \; | \; \Re (s) > \sigma\}\) and
let \(Z_{\eta}\) be the set of points whose distance to the zeros of
\(\det \Delta_0\) is at most \(\eta\).}

\emph{For any \(\sigma \in \mathbb{R}\) and \(\epsilon > 0\), there is a
\(\eta>0\) such that:}

\begin{enumerate}
\def\labelenumi{\roman{enumi}.}
\item
  \emph{any connected component \(\Lambda\) of the set \(Z_{\eta}\) is
  bounded,}
\item
  \emph{\(\Lambda\) is a subset of \(H_{\sigma - \epsilon}\) whenever
  \(\Lambda \cap H_{\sigma} \neq \varnothing\),}
\item
  \emph{\(|\det \Delta_0|\) has a positive lower bound on
  \(H_{\sigma-\epsilon} \setminus Z_{\eta}\).}
\end{enumerate}

The proof (from (Boisgérault \protect\hyperlink{ref-Boi13}{2013})) is
given in the appendix for the sake of completeness.

\textbf{Lemma -- Characteristic Function Zeros.} \label{ateam} \emph{Let
\(\sigma \in \mathbb{R}\). If the function \(\det \Delta_0\) has an
infinite number of zeros on \(H_{\sigma}\), the function
\(\det \Delta^{\dagger}\) has an infinite number of zeros on
\(H_{\sigma - \epsilon}\) for any \(\epsilon>0\).}

\textbf{Proof.} Suppose that \(\det \Delta_0\) has an infinite number of
zeros in \(H_{\sigma}\) and let \(\eta>0\) be as in the zero clusters
lemma. The zeros of \(\det \Delta_0\) are isolated, thus the collection
of components \(\Lambda\) of \(Z_{\eta}\) such that
\(\Lambda \cap H_{\sigma} \neq \varnothing\), which is locally finite by
construction, is also infinite: for any compact set \(K\), there is a
\(\Lambda\) such that \(\Lambda \cap K = \varnothing\).

Since \(\det \Delta^{\dagger}(s)\) has a quasi-polynomial structure with
leading coefficient \(\det \Delta_0(s) \, s^n\) and since
\(\det \Delta_0(s)\) has a positive lower bound on
\(H_{\sigma-\epsilon} \setminus Z_{\eta}\), there is a compact set \(K\)
such that
\(|\det \Delta^{\dagger}(s) - s^n \det \Delta_0(s)| < |s^n \det \Delta_0(s)|\)
whenever \(s \in H_{\sigma-\epsilon} \setminus Z_{\eta}\) and
\(s \not \in K\). Rouché's theorem is applicable to any of the
components \(\Lambda\) in the complement of \(K\), each of which
contains at least one zero of \(\det \Delta_0\). Thus,
\(\det \Delta^{\dagger}\) has an infinite number of zeros in
\(H_{\sigma-\epsilon}\). \hfill \(\blacksquare\)

We can associate to the initial-value problem \eqref{DDAE}, or
\eqref{DDAEX} which is equivalent, a one-parameter family of operators
\((\exp(\mathcal{A}t))_{t \geq 0}\) defined by \[
   (x(t^+), x_t, y_t) = \exp (\mathcal{A} t) (\phi, \chi, \psi)
   \; \mbox{ for } \; t \geq 0 , \; (\phi,\chi,\psi) \in X.
  \] It is a strongly continous semigroup on the product space \(X\);
its infinitesimal generator \(\mathcal{A}\) is defined by \[
\mathcal{A}(\phi,\chi,\psi) 
=
(
  A\chi + B \psi, \; \dot{\chi}, \; \dot{\psi}
)
\] on the domain \[
  \{
    (\phi,\chi,\psi) \in \mathbb{C}^n \times W^{1,2}([-r,0],\mathbb{C}^{n+m}) 
    \; | \; 
    \chi(0) = \phi, \; \psi(0) = E\chi + F\psi
  \}.
  \]

The resolvent operator \((sI - \mathcal{A})^{-1}\) exists if and only if
\(\Delta^{\dagger}(s)\) is non-singular and moreover, for any real
number \(\sigma\) there are constants \(\kappa_{\sigma}\) and
\(\lambda_{\sigma}\) such that

\begin{equation} \label{bbb}
  \|(sI-\mathcal{A})^{-1}\| 
  \leq 
  \kappa_{\sigma} \|\Delta^{\dagger}(s)^{-1}\| + \lambda_{\sigma}
  \end{equation}

if \(\Re s \geq \sigma\) and \(\det \Delta^{\dagger}(s) \neq 0\) (see
Salamon (\protect\hyperlink{ref-Sal84}{1984})).

We may now prove the main result:

\textbf{Proof.} Let \(s(\mathcal{A})\) be the spectral bound of
\(\mathcal{A}\). We show that for any \(\sigma > s(\mathcal{A})\),
\(\|\Delta^{\dagger}(s)^{-1}\|\) is bounded on
\(\overline{H}_{\sigma}\). Thus, by inequality \eqref{bbb},
\(\|(sI-\mathcal{A})^{-1}\|\) is similarly bounded and the
Gearhart-Prüss theorem proves the result.

The quasi-polynomial structure of the adjugate matrix provides for some
\(\kappa \geq 0\) \[
  \|\mbox{\rm adj} \, \Delta^{\dagger}(s)\| \leq \kappa (1 + |s|^n)
  \] on \(\overline{H}_{\sigma}\). Since the resolvent operator is
defined on \({H}_{\sigma - \epsilon}\) for any \(\epsilon>0\) such that
\(s(\mathcal{A}) < \sigma - \epsilon\), \(\det \Delta^{\dagger}\) has no
zero in \({H}_{\sigma - \epsilon}\). Thus, by the characteristic
function zeros lemma, \(\det \Delta_0\) has at most a finite number of
zeros on \({H}_{\sigma - \epsilon / 2}\) and thus on
\(\overline{H}_{\sigma}\). By the zero clusters lemma, on
\(\overline{H}_{\sigma}\) and away from these zeros, \(|\det \Delta_0|\)
has a positive lower bound \(\kappa'\). It follows from the
quasipolynomial structure of \(\det \Delta^{\dagger}\) that \[
|\det \Delta^{\dagger} (s)| \geq \frac{\kappa'}{2} (1 + |s|^n)
\] on \(\overline{H}_{\sigma}\) except on some compact set \(K\); by
continuity of \(\det \Delta^{\dagger}\), this estimate still holds on
all of \(\overline{H}_{\sigma}\) with a possibly smaller \(\kappa'\).
Finally, for any \(s\in\overline{H}_{\sigma}\),
\(\|\Delta^{\dagger} (s)^{-1}\| = {\|\mbox{\rm adj} \, \Delta^{\dagger}(s)\|} / {|\det \Delta^{\dagger}(s)|} \leq 2\kappa / \kappa'.\)
\hfill \(\blacksquare\)

\subsection{Stability of the FSA
System}\label{stability-of-the-fsa-system}

Since the kernels of FSA system \eqref{FSA} are given by
\eqref{FSA-ker-1} and \eqref{FSA-ker-2}, its characteristic matrix
satisfies \[
  \renewcommand\arraystretch{1.3}
  \Delta(s)
  =
  \left[
  \begin{array}{c|c}
  sI_n - E & FG \\
  \hline
  -e^{(sI_n - E)T} & I_n + P(s) FG
  \end{array}
  \right]
  \] where \(P(s)\) is the analytic extension to \(\mathbb{C}\) of the
meromorphic function \[
  s \mapsto [sI_n - E]^{-1}(I_n - e^{(sI_n-E)T}).
  \] A straightforward computation provides \[
  \Delta(s)
  \left[   
  \begin{array}{c|c}
  I_n & 0 \\
  \hline
  I_n & I_n
  \end{array}
  \right]
  =
  \left[
  \begin{array}{c|c}
  sI_n - E + FG & FG \\
  \hline
  P(s)(s I_n - E + FG) & I_n +P(s) FG
  \end{array}
  \right]
  \] and the right-hand side of this equation can be factored into \[
  \left[
  \begin{array}{c|c}
  I_n & 0 \\
  \hline
  P(s) & I_n
  \end{array}
  \right]
  \left[
  \begin{array}{c|c}
  I_n & FG \\
  \hline
  0 & I_n
  \end{array}
  \right]
  \left[
  \begin{array}{c|c}
  sI_n - E + FG & 0 \\
  \hline
  0 & I_n
  \end{array}
  \right].
  \] Thus, the characteristic function satisfies \[
  \det \Delta (s) = \det (sI_n - E + FG)
  \] and it roots are the eigenvalues of \(E - FG\). When the delay-free
system \(\dot{x}(t) = E x(t) + F u(t)\) is controllable,
finite-dimensional control theory provides a gain matrix \(G\) that
assigns these eigenvalues to \(n\) arbitrary locations in the complex
plane, and thus we may achieve a ``finite spectrum assignment'' of the
closed-loop system. In particular, it's possible to select a gain matrix
\(G\) which exponentially stabilizes system \eqref{FSA}.

\section*{Appendix}\label{appendix}
\addcontentsline{toc}{section}{Appendix}

We provide in this appendix a proof of the zero clusters lemma. First we
derive some properties of \(\det \Delta_0\) by expressing it as the
Laplace transform of a complex measure \(\mu\). Let \(\Sigma_m\) be the
set of permutations of \(\{1,...,m\}\) and
\[{\rm det}_{\ast} \, M = \sum_{\sigma \in \Sigma_m} 
\mathrm{sgn}(\sigma) M_{1,\sigma(1)} \ast \hdots \ast M_{{m},\sigma(m)}.\]
As \(\Delta_0(s) = I_m - \mathcal{L}F^{\ast}(s)\),
\(\det \Delta_0 = \mathcal{L} \mu\) where
\(\mu = {\rm det}_{\ast} (\delta_0 I_m - F^{\ast}).\) The complex
measure \(\mu\) is a sum of convolution products of \(m\) complex
measures supported on \([0,r]\), hence it is supported on \([0,mr]\).
Consequently, \(\det \Delta_0\) is an entire function that satisfies the
inequality

\begin{equation} \label{exp-bound}
  |\det \Delta_0(s)| \leq |\mu|([0,mr]) \max(1, \exp(-\Re(s) m r)).
  \end{equation}

Since \(F \{0\} = 0\), we also have \(\mu \{0\} = 1\), which yields

\begin{equation} \label{limit-one}
\lim_{\Re s \to  +\infty} \det \Delta_0(s) = 1.
\end{equation}

We may now proceed to the proof:

\textbf{Proof (Zero Clusters Lemma).} The function
\(z \mapsto \det \Delta_0(iz)\) meets the assumptions of theorem VIII in
(Levinson \protect\hyperlink{ref-Lev40}{1940}) and is of exponential
type \(mr\). Thus, the number of distinct zeros \(N(\rho)\) of
\(\det \Delta_0\) whose modulus is less than \(\rho\) is such that
\(\limsup_{\rho \to +\infty} {N(\rho)} / {\rho} \leq {2mr} /{\pi}.\) If
the component of \(Z_{\eta}\) that contains a zero \(s\) is unbounded,
there are at least \(n+1\) zeros in the closed disk centered at \(s\) of
radius \(2\eta n\) and thus
\(\limsup_{\rho \to +\infty} {N(\rho)}/{\rho} \geq 1/2\eta\).
Consequently, if \(\eta < \pi / 4mr\), every connected component of
\(Z_{\eta}\) is bounded.

The proofs of statements ii.~and iii.~use a similar argument. In each
case we consider a sequence \(s_n\) of complex numbers of bounded real
part and the functions \(f_n(s) = \det \Delta_0(s + i \Im s_n)\). Since
inequality (\ref{exp-bound}) holds, these functions are locally
uniformly bounded, hence there is a subsequence of \(s_n\) which
converges to some real number \(x\) and, by Montel's theorem, a
corresponding subsequence of \(f_n\) that converges locally uniformly to
an entire function \(f_{\infty}\). By (\ref{limit-one}), \(f_{\infty}\)
is not identically zero, thus by Hurwitz's theorem, if \(m\) is the
multiplicity of \(x\) if \(f_{\infty}(x)=0\), or \(0\) otherwise, for
any sufficiently small \(\delta>0\), \(\det \Delta_0\) has exactly \(m\)
zeros in the open disk \(B(s_n, \delta)\) for an infinite number of
values of \(n\).

To prove ii., we assume the existence of a sequence \(\Lambda_n\) of
connected components of \(Z_{\eta_n}\) where \(\eta_n \to 0\) and such
that \(\Lambda_n \cap H_{\sigma} \neq \varnothing\) but
\(\Lambda_n \not \subset H_{\sigma - \varepsilon}.\) The \(\Lambda_n\)
are eventually bounded and there is a tuple of zeros in \(\Lambda_n\)
with a first element such that \(\Re s + \eta_n > \sigma\), a last
element such that \(\Re s - \eta_n \leq \sigma - \epsilon\) and a
distance between consecutive points that is at most \(2\eta_n\). Let
\(s_n\) to be the first element of this tuple; the real part of this
sequence is bounded by \eqref{limit-one}. But for any \(\delta > 0\),
the number of zeros in \(B(s_n, \delta)\) converges to \(+\infty\), a
contradiction with the previous paragraph.

To prove iii., we assume the existence of a sequence \(s_n\) in
\(H_{\sigma-\epsilon} \setminus Z_{\eta}\) such that
\(\det \Delta_0(s_n) \to 0\) when \(n \to +\infty\); the real part of
this sequence is bounded and we may select \(s_n\) such that
\(\Re s_n \to x\). Now,
\(f_{\infty}(x) = \lim_{n\to+\infty} f_n(\Re s_n) = 0\), thus for some
\(\delta < \eta\) and for some value of \(n\), there is at least one
zero of \(\det \Delta_0\) in \(B(s_n,\delta)\), which is a
contradiction. \hfill \(\blacksquare\)

\section*{References}\label{references}
\addcontentsline{toc}{section}{References}

\hypertarget{refs}{}
\hypertarget{ref-BPDM06}{}
Bensoussan, Alain, Giuseppe Da Prato, Michel C. Delfour, and Sanjoy K.
Mitter. 2006. \emph{Representation and Control of Infinite Dimensional
Systems (Systems \& Control: Foundations \& Applications)}. Birkhauser.

\hypertarget{ref-Boi13}{}
Boisgérault, Sébastien. 2013. ``Growth bound of delay-differential
algebraic equations.'' \emph{C. R., Math., Acad. Sci. Paris} 351
(15-16). Elsevier (Elsevier Masson), Issy-les-Moulineaux; Académie des
Sciences, Paris: 645--48.
doi:\href{https://doi.org/10.1016/j.crma.2013.08.001}{10.1016/j.crma.2013.08.001}.

\hypertarget{ref-DK87}{}
Delfour, M.C., and J. Karrakchou. 1987. ``State space theory of linear
time invariant systems with delays in state, control, and observation
variables. I, II.'' \emph{Journal of Mathematical Analysis and
Applications}.

\hypertarget{ref-FD88}{}
Fell, J.M.G., and R.S. Doran. 1988. \emph{Representations of *-algebras,
locally compact groups, and Banach *- algebraic bundles. Vol. 1: Basic
representation theory of groups and algebras.} Boston, MA etc.: Academic
Press, Inc.

\hypertarget{ref-H77}{}
Hale, J. K. 1977. \emph{Theory of Functional--differential Equations}.
Springer--Verlag, Berlin--Heidelberg--New York.

\hypertarget{ref-Lev40}{}
Levinson, Norman. 1940. \emph{Gap and density theorems.} American
Mathematical Society (AMS). Colloquium Publications. 26. New York:
American Mathematical Society (AMS). VIII, 246 p.

\hypertarget{ref-MO79}{}
Manitius, Andrzej Z., and Andrzej W. Olbrot. 1979. ``Finite spectrum
assignment problem for systems with delays.'' \emph{IEEE Trans. Autom.
Control} 24: 541--53.
doi:\href{https://doi.org/10.1109/TAC.1979.1102124}{10.1109/TAC.1979.1102124}.

\hypertarget{ref-MR03}{}
Mirkin, Leonid, and Natalya Raskin. 2003. ``Every stabilizing dead-time
controller has an observer-predictor-based structure.''
\emph{Automatica} 39 (10): 1747--54.
doi:\href{https://doi.org/10.1016/S0005-1098(03)00182-1}{10.1016/S0005-1098(03)00182-1}.

\hypertarget{ref-Sal84}{}
Salamon, D. 1984. \emph{Control and observation of neutral systems.}
Research Notes in Mathematics, 91, Boston-London-Melbourne: Pitman
Advanced Publishing Program. 207 p.

\end{document}